\documentclass[10pt]{article}
\usepackage[cp1251]{inputenc}
\usepackage[ukrainian, russian, english]{babel}
\usepackage{ amsfonts, amssymb}

\begin{document}
 \selectlanguage{ukrainian} \thispagestyle{empty}
 \pagestyle{myheadings}                 

УДК 517.51 \vskip 3mm

\bf{А.\,Ф.~Конограй} \rm
\footnote{Работа выполнена при частичной поддержке Государственного фонда фундаментальных исследований Украины (проект \mbox{№ GP/Ф36/068)}} \rm (Институт математики НАН Украины, Киев)

\vskip 7mm

\bf{\centerline{Оценки аппроксимативных характеристик классов
\boldmath{$B^{\Omega}_{p,\theta}$}} \centerline{периодических
функций многих переменных} \centerline{с заданной мажорантой
смешанных модулей непрерывности}}\rm

\vskip 7mm

\emph{В работе получены точные по порядку оценки ортопроекционных
поперечников классов $B^{\Omega}_{p,\theta}$ периодических функций
многих переменных, мажоранта смешанных модулей непрерывности которых
содержит как степенные, так и логарифмические множители.}

\vskip 7mm

\bf{1. Введение.} \rm В настоящей работе продолжаются исследования
аппроксимативных характеристик классов $B^{\Omega}_{p,\theta}$
периодических функций многих переменных с заданной функцией
$\Omega(t)$ специального вида, которые изучались в работах [\ref{1},
\ref{2}]. Более конкретно об этом будет идти речь в соответствующих
частях работы, а сначала приведем необходимые обозначения и
определения.

Пусть $\mathbb{R}^{d}, d\geq1,$ обозначает $d$-мерное пространство с
элементами

$$
{x=(x_1,...,x_d)}, \ \ (x,y)=x_{1}y_{1}+...+x_{d}y_{d}
$$
и
$$
L_p(\pi_d), \ \ {\pi_{d}=\prod\limits_{j=1}^{d}[0;2\pi]},
$$
--- пространство $2\pi$-периодических по каждой переменной функций
$f$, для которых
$$
\|f\|_p:=
\bigg(\left(2\pi\right)^{-d}\int\limits_{\pi_{d}}|f(x)|^{p}dx \bigg)
^{\frac{1}{p}}<\infty, \ \ \ 1\leq p<\infty,
$$
$$
\|f\|_{\infty}:=\mathop {\rm ess \sup}\limits_{x\in \pi_d}
|f(x)|<\infty.
$$

Далее будем предполагать, что для функций $f\in L_p(\pi_d)$
выполнено дополнительное условие
 $$ \int\limits_{0}^{2\pi}f(x)dx_j=0 \ ,\ \ j=\overline{1,d} .$$

Для\ \ $f\in L_p(\pi_d), 1\leq p\leq\infty$, и  $t=(t_1,...,t_d)$,
$t_j\geq0$, $j=\overline{1,d}$, рассмотрим смешанный модуль
непрерывности порядка \ $l$
$$
\Omega_l(f,t)_p=\sup_{\scriptstyle |h_j|\le t_j\atop \scriptstyle
j=\overline{1,d}}||\Delta_h^lf(\cdot)||_p,
$$
где $l\in \mathbb{N}, \
\Delta_h^lf(x)=\Delta_{h_1}^l\dots\Delta_{h_d}^lf(x)$
${=\Delta_{h_d}^l(\dots(\Delta_{h_1}^lf(x)))}$ --- смешанная
разность порядка $l$ с шагом  $h_j$ по переменной $x_j$ и
$$
\Delta_{h_j}^lf(x)=\sum\limits_{n=0}^l(-1)^{l-n}C_l^n
f(x_1,\dots,x_{j-1},x_j+nh_{j},x_{j+1},\dots,x_d).
$$

Пусть $\Omega(t)=\Omega(t_1,\dots,t_d)$ --- заданная функция типа
смешанного модуля непрерывности порядка $l$, которая удовлетворяет
следующим условиям:

1) $\Omega(t)>0,$ $t_j>0,$ $j=\overline{1,d};$ $\Omega(t)=0,$
$\prod \limits_{j=1}^dt_j=0;$

2) $\Omega(t)$ не убывает по каждой переменной;

3) $\Omega(m_1t_1,\dots,m_dt_d)\leq\bigg(\prod
\limits_{j=1}^dm_j\bigg)^l\Omega(t),$  $m_j\in \mathbb{N},$ \ \
$j=\overline{1,d};$

4) $\Omega(t)$\ \ непрерывна при\ \ $t_j\geq 0,$ \ \
$j=\overline{1,d}$ .

Будем говорить, что $\Omega(t)$ удовлетворяет также условия $(S)$ и
$(S_l),$ которые называют условиями Бари--Стечкина [\ref{3}]. Это
означает следующее.

Функция одной переменной $\varphi(\tau)\geq 0$ удовлетворяет условию
$(S),$ если $\varphi(\tau)/\tau^{\alpha}$ почти возрастает при
некотором $\alpha>0,$ т.е. существует независимая от \ $\tau_1$\ \
и\ \ $\tau_2 $ постоянная $C_1>0,$ такая, что
 $$
 \frac{\varphi(\tau_1)}{\tau_1^{\alpha}} \leq
C_1\frac{\varphi(\tau_2)}{\tau_2^{\alpha}} \ , \ \ \ \ 0<\tau_1 \le
\tau_2 \le 1 .
 $$

Функция \  $\varphi(\tau)\geq 0$ удовлетворяет условию $(S_l),$ если
при некотором ${0<\gamma<l}$, $\varphi(\tau)/\tau^{\gamma}$ почти
убывает, т.е. существует независимая от \ $\tau_1$\ \ и\ \ $\tau_2 $
постоянная $C_2>0 ,$ такая, что
$$
\frac{\varphi(\tau_1)}{\tau_1^{\gamma}} \geq
C_2\frac{\varphi(\tau_2)}{\tau_2^{\gamma}} \ , \ \ \ \ 0<\tau_1 \le
\tau_2 \le 1 .
$$

Будем считать, что $\Omega(t)$ удовлетворяет условиям $(S)$ и
$(S_l)$, если $\Omega(t)$ удовлетворяет эти условия по каждой
переменной \ $t_j$ при фиксированных $t_i ,$ \ $i\neq j \ .$

Итак, пусть $1\leq p\leq\infty$, ${1\le\theta\le\infty}$, а
$\Omega(t)$ --- заданная функция типа смешанного модуля
непрерывности порядка $l$, которая удовлетворяет условиям 1~--~4.
Тогда (см.~[\ref{4}])
 $$
B_{p,\theta}^{\Omega}:=\!\bigg\{f\in  L_p(\pi_d):
\|f\|_{B_{p,\theta}^{\Omega}}= \bigg (\int \limits_{\pi_d} \bigg
(\frac{\Omega_l(f,t)_p}{\Omega(t)} \bigg)^{\theta} \prod
\limits_{j=1}^d\frac{dt_j}{t_j} \bigg)^{\frac{1}{\theta}}\leq
1\bigg\},
 $$
при $1\le\theta<\infty$ и
 $$
B_{p,\infty}^{\Omega}:=\bigg\{f\in L_p(\pi_d): \
\|f\|_{B_{p,\infty}^{\Omega}}=\sup\limits_{t>0}
\frac{\Omega_l(f,t)_p}{\Omega(t)}\leq 1\bigg\},
 $$
(запись $t>0$ для $t=(t_1,...,t_d)$ равносильна $t_{j}>0,
j=\overline{1,d}$).

Заметим, что при \ $\theta=\infty$\ \ классы\ \
$B_{p,\theta}^{\Omega}$ совпадают с классами \ $H_p^{\Omega}$,
которые были рассмотрены Н.Н. Пустовойтовым в [\ref{5}].

В последующих рассуждениях нам будет удобно пользоваться
эквивалентным (с точностью до абсолютных постоянных) определением
классов $B_{p,\theta}^{\Omega}$.

Каждому вектору $s=(s_{1},...,s_{d}), \ s_{j}\in \mathbb{N}, \
j=\overline{1,d}$, поставим в соответствие множество
 $$
\rho(s)=\left\{k=(k_1,...,k_d): 2^{s_j-1}\le|k_j|<2^{s_j},
k_{j}\in\mathbb{Z}\setminus\{0\}, j=\overline{1,d}\right\}
 $$
и для $f\in L_p(\pi_d)$ обозначим
 $$
\delta_s(f,x)=\sum\limits_{k\in \rho(s)}\widehat{f}(k)e^{i(k,x)},
 $$
где
$$\widehat{f}(k)=(2\pi)^{-d}\int\limits_{\pi_d}f(t)e^{-i(k,t)}dt$$
--- коэффициенты Фурье функции $f$.

Итак, пусть $1< p<\infty$, \ $1\le\theta\le\infty$\ и $\Omega(t)$
--- заданная функция типа смешанного модуля непрерывности порядка $l$,
которая удовлетворяет условиям 1 --- 4, $(S)$ и $(S_l),$ тогда с
точностью до абсолютных постоянных классы $B_{p,\theta}^{\Omega}$
можно определить следующим образом (см. [\ref{4}]):
\begin{equation}
B_{p,\theta}^{\Omega}:=\bigg\{f\in L_p(\pi_d):
\|f\|_{B_{p,\theta}^{\Omega}}=\bigg(\sum\limits_s
\Omega^{-\theta}(2^{-s}) \|\delta_s(f,\cdot)\|_p^\theta\
\bigg)^{\frac{1}{\theta}}\leq1\bigg\}, \label{1.1}
\end{equation}
при $1\le\theta<\infty$ и
\begin{equation}
B_{p,\infty}^{\Omega}:=\bigg\{f\in L_p(\pi_d):
\|f\|_{B_{p,\infty}^{\Omega}}=\sup\limits_s
\frac{\|\delta_s(f,\cdot)\|_p}{\Omega(2^{-s})} \leq 1\bigg\},
\label{1.2}
\end{equation}
здесь и в дальнейшем $\Omega(2^{-s})=\Omega(2^{-s_1},...,2^{-s_d}),$
$s_j\in \mathbb{N} ,$ $j=\overline{1,d}.$

Отметим, что при \mbox{$\Omega(t)=\prod\limits_{j=1}^dt_j^{r_j}, \
0<r_j<l$}, классы $B_{p,\theta}^{\Omega}$ есть аналоги известных
классов Бесова $B_{p,\theta}^{r}$ (см., например, [\ref{6}]).

Приведенное определение классов $B^{\Omega}_{p,\theta}$ можно
распространить и на крайние значения $p=1, \infty$, видоизменив в
(\ref{1.1}) и (\ref{1.2}) "блоки"\ $\delta_s(f,x)$.

Пусть $V_n(t), n\in \mathbb{N},$ обозначает ядро Валле Пуссена
порядка $2n-1$:
 $$
V_n(t)=1+2\sum_{k=1}^{n}\cos
kt+2\sum_{k=n+1}^{2n-1}\bigg(1-\frac{k-n}{n}\bigg)\cos kt .
 $$

Сопоставим каждому вектору $s=(s_1,\ldots,s_d),$ $s_j\in\mathbb{N},$
$j=\overline{1,d},$ полином
$$
A_s(x)=\prod_{j=1}^d\bigg(V_{2^{s_j}}(x_j)-V_{2^{s_j-1}}(x_j)\bigg)
$$
и для $f\in L_p(\pi_d),$ $1\leq p\leq\infty$, через $A_s(f,x)$
обозначим свертку
 $$
A_s(f,x)=f\ast A_s.
 $$

В принятых обозначениях (с точностью до абсолютных постоянных)
классы $B_{p,\theta}^{\Omega},$ $1\leq p\leq\infty$ определяются
следующим образом (см., соответственно [\ref{7}] и~[\ref{5}]):
\begin{equation}
B_{p,\theta}^{\Omega}:=\bigg\{f\in L_p(\pi_d):
 \|f\|_{B_{p,\theta}^{\Omega}}=
 \bigg(\sum\limits_{s} \Omega^{-\theta}(2^{-s})
\left\|A_s(f,\cdot)\right\|_p^\theta
\bigg)^{\frac{1}{\theta}}\leq1\bigg\}, \label{1.3}
\end{equation}
при $1\le\theta<\infty$ и
 \begin{equation} B_{p,\infty}^{\Omega}:=\bigg\{f\in
L_p(\pi_d): \|f\|_{B_{p,\infty}^{\Omega}}=\sup_{s}
\frac{\left\|A_s(f,\cdot)\right\|_p}{\Omega(2^{-s})}\leq 1\bigg\} .
\label{1.4}
 \end{equation}

Отметим, что при $1<p<\infty$ определения норм функций с классов
$B_{p,\theta}^{\Omega}$ (\ref{1.1}) и (\ref{1.2}) эквивалентны к
определениям норм (\ref{1.3}) и (\ref{1.4}) соответственно.

Ниже мы будем рассматривать классы $B_{p,\theta}^{\Omega}$, которые
определяются функцией $\Omega(t)$ специального вида, а именно
\begin {equation}
\Omega(t)=\Omega(t_{1},...,t_{d})= \left\{\begin{array}{cc}
\displaystyle \prod\limits_{j=1}^{d}
\frac{t_{j}^{r}}{\big(\log\frac{1}{t_{j}}\big)_{+}^{b_{j}}}, &
\mbox{если} \  t_{j}>0, \  j=\overline{1,d};
\\ \displaystyle 0,\ & \mbox{если} \ \prod\limits_{j=1}^{d}t_{j}=0.
\end{array} \right.
\label{1.5}
 \end {equation}

Здесь
$$
\bigg(\log\frac{1}{t_{j}}\bigg)_{+}=\max\bigg\{1,
\log\frac{1}{t_{j}}\bigg\},
$$
причем логарифм, как и всюду ниже, берется по основанию 2.

Также полагаем, что $b_{j}<r, j=\overline{1,d}$, и $0<r<l$, а значит
функция $\Omega(t)$, заданная формулой (\ref{1.5}), будет
удовлетворять сформулированным выше условиям 1--4, $(S)$ и $(S_l)$.

В работе получены точные по порядку оценки ортопроекционных
поперечников классов $B_{p,\theta}^{\Omega}$ в пространстве $L_{q}$
при некоторых соотношениях между параметрами $p$ и $q$, а также
установлены порядки приближения этих же классов функций в
пространстве $L_{q}$ с помощью линейных операторов, которые
подчинены определенным условиям. Соответствующие аппроксимативные
характеристики связаны с ортопроекционными поперечниками классов
$B_{p,\theta}^{\Omega}$, о чем более конкретно будет идти речь ниже.

Для формулировки и доказательства полученных результатов нам
понадобятся соответствующие обозначения и определения.

Пусть $\{u_{i}(x)\}_{i=1}^{M}$ --- ортонормированная система функций
$u_{i}\in L_{\infty}(\pi_{d})$. Каждой функции $f\in L_{q}(\pi_{d}),
\ 1\leq q\leq\infty$, поставим в соответствие приближающий аппарат
вида $\sum\limits_{i=1}^{M}(f,u_{i})u_{i}(x)$, т.е. ортогональную
проекцию функции $f$ на подпространство, порожденное системой
функций $\{u_{i}(x)\}_{i=1}^{M}$. Если $F\subset L_{q}(\pi_{d})$ ---
класс функций, то величина
\begin{equation}
d_M^{\bot}(F,L_q)=\inf_{\{u_i(x)\}^M_{i=1}}\sup_{f\in
F}\bigg\|f(\cdot)-\sum_{i=1}^M(f,u_i)u_i(\cdot)\bigg\|_q,
\label{1.6}
\end{equation}
называется \ ортопроекционным \ поперечником \ (Фурье-поперечником)
\ класса $F$ в пространстве $L_{q}(\pi_{d})$. Поперечник
$d_M^{\bot}(F,L_q)$ введен В.Н. Темляковым в работе [\ref{8}], где
также рассмотрена величина
\begin{equation}
d_M^{B}(F,L_q)=\inf_{G\in L_{M}(B)_{q}}\sup_{f\in F\cap
D(G)}\|f(\cdot)-Gf(\cdot)\|_q. \label{1.7}
\end{equation}
Здесь через $L_{M}(B)_{q}$ обозначено множество линейных операторов,
удовлетворяющих условиям:

а) область определения $D(G)$ этих операторов содержит все
тригонометрические полиномы, а их область значений содержится в
подпространстве размерности $M$ пространства $L_q(\pi_d);$

б) существует число $B\geq1$ такое, что для всех векторов
$k=(k_{1},...,k_{d})$, выполнено неравенство
$$
\|Ge^{i(k,\cdot)}\|_{2}\leq B.
$$

Отметим, что к $L_{M}(1)_{2}$  принадлежат операторы ортогонального
проектирования на пространства размерности $M$, а также операторы,
которые задаются на ортонормированной системе функций с помощью
мультипликатора, определяющегося последовательностью
$\{\lambda_{l}\}$ такой, что $|\lambda_{l}|\leq1$ для всех $l$.
Легко видеть, что согласно определений
$$
 d_M^{B}(F,L_q)\leq d_M^{\bot}(F,L_q).
$$

Следовательно оценки \ величин \ $d_M^{B}(F,L_q)$ могут \ служить
оценками снизу для ортопроекционных \ поперечников \
$d_M^{\bot}(F,L_q)$ и наоборот -- оценки \ поперечников
$d_M^{\bot}(F,L_q)$ можно использовать для оценок сверху величин
$d_M^{B}(F,L_q)$. Это обстоятельство будет нами использоваться при
доказательстве соответствующих утверждений.

Отметим, что величины (\ref{1.6}) и (\ref{1.7}) некоторых классов
функций исследовались  в работах [\ref{9} -- \ref{15}] (см. также
монографии [\ref{16}, \ref{17}]), где можно ознакомиться и с
соответствующей библиографией.

В процессе доказательства полученных результатов используются и
развиваются методы, которые применялись при исследовании
рассматриваемых аппроксимативных характеристик в работах [\ref{8} --
\ref{17}]. Более подробно об этом будет говориться в соответствующих
комментариях.

\bf{2. Вспомогательные утверждения.}\rm

Напомним несколько известных утверждений, которые нами будут
систематически использоваться.

{\bf Теорема А} (Литтлвуда--Пэли~[\ref{18}]). {\it Пусть задано
$p\in(1, \infty)$. Существуют положительные константы $C_{3}(p)$ и
$C_{4}(p)$ такие, что для каждой функции ${f\in L_{p}(\pi_{d})}$
справедлива оценка
 $$
 C_{3}(p)\|f\|_{p}\leq \Big\| \Big(\sum\limits_{s}|\delta_{s}(f,\cdot)|^{2} \Big)
 ^{\frac{1}{2}} \Big\|_{p}\leq C_{4}(p)\|f\|_{p}.
 $$}

Непосредственно из теоремы А легко получить следующее

{\bf Следствие.} {\it
 Для $f\in L_{p}(\pi_{d})$ при $1<p<\infty$ имеет место следующее соотношение
 $$
 \|f\|_{p}\ll \Big(\sum\limits_{s}\|\delta_{s}(f,\cdot)\|^{p_0}_p \Big)
 ^{\frac{1}{p_0}},
 $$
где $p_0=\min\{p;2\}$.}

Здесь и далее, для положительных функций $\mu_1(N)$ и $\mu_2(N)$
запись $\mu_1\ll\mu_2$ означает, что  существует постоянная $C>0$
такая, что $\mu_1(N)\leq C\mu_2(N)$. Соотношение $\mu_1\asymp \mu_2$
равносильно тому, что $\mu_1\ll\mu_2$ и $\mu_1\gg\mu_2$. Отметим,
что все постоянные $C_i, i=1,2,\ldots,$ которые будут встречаться в
работе могут зависеть только от тех параметров, которые содержатся
 в определениях классов, метрики и размерности пространства
$\mathbb{R}^d$.

{\bf Теорема Б \rm [\ref{19}].} {\it Пусть ${n=(n_{1},...,n_{d})}, \
n_{j}$ --- целые неотрицательные числа, $j=\overline{1,d}$, и
 $$
t(x)=\sum\limits_{|k_{j}|\leq n_{j}} c_{k}e^{i(k,x)}.
 $$
Тогда при $1 \leq q < p\leq\infty$ имеет место неравенство
 $$
\|t\|_{p}\leq 2^{d} \prod \limits_{j=1}^{d}
n_{j}^{\frac{1}{q}-\frac{1}{p}}\|t\|_{q}.
 $$}
Это неравенство установлено С.М. Никольским и названо “неравенством
разных метрик”. В одномерном случае и при  $p=\infty$
соответствующее неравенство доказал Джексон [\ref{20}].

{\bf Лемма A \rm[\ref{10}].} {\it Пусть линейный оператор $A$ задан
равенством
$$
Ae^{i(k,x)}=\sum\limits_{m=1}^{L}a_m^k \psi_{m}(x),
$$
где $\big\{\psi_m(x)\big\}_{m=1}^{L}$ --- набор функций, для которых
$$
\|\psi_m(\cdot)\|_2\leq1, \ \ m=1,...,L.
$$
Тогда для любого тригонометрического полинома $t$ верно неравенство}
$$
\min\limits_{y=x} \mbox{Re} \ At(x-y)\leq
\bigg(L\sum\limits_{m=1}^{L}\sum\limits_k|a_m^k\widehat{t}(k)|^2
\bigg)^{\frac{1}{2}}.
$$

Для натурального $N$ положим
 $$
\chi(N)=\bigg\{s=(s_{1},..., s_{d}): \ s_{j}\in\mathbb{N}, \
j=\overline{1,d}, \ \Omega(2^{-s})\geq\frac{1}{N}\bigg\}.
 $$

С учетом определения (\ref{1.5}) множество $\chi(N)$ можно
определить следующим образом:
 $$
\chi(N)=\bigg\{s: \ s_{j}\in\mathbb{N}, \ j=\overline{1,d}, \
\prod\limits_{j=1}^{d}2^{rs_j}s_j^{b_j}\leq N\bigg\}.
 $$

Далее, пусть
 $$
\chi^{\perp}(N)=\mathbb{N}^{d}\setminus\chi(N)
 $$

В качестве приближающих подпространств возьмем подпространства
тригонометрических полиномов, спектр которых содержится во
множествах
 $$
Q(N)=\bigcup\limits_{s\in\chi(N)}\rho(s).
 $$

Также нам понадобятся множества
 $$
\Theta(N)=\bigg\{s=(s_{1},..., s_{d}): \ s_{j}\in\mathbb{N}, \
j=\overline{1,d}, \
\frac{1}{2^{l}N}\leq\Omega(2^{-s})<\frac{1}{N}\bigg\}.
 $$

В принятых обозначениях имеют место следующие утверждения, которые
получены Н.Н. Пустовойтовым.

{\bf Лемма Б \rm [\ref{14}].} {\it Количество элементов множества
$Q(N)$ по порядку равно:}
 $$
|Q(N)|\asymp N^{\frac{1}{r}}\big(\log
N\big)^{-\frac{b_{1}}{r}-...-\frac{b_{d}}{r}+d-1}.
 $$

{\bf Лемма В \rm [\ref{21}].} {\it Для количества элементов
множества $\Theta(N)$ имеет место соотношение}
$$
|\Theta(N)|\asymp (\log N)^{d-1}.
$$

{\bf Лемма Г \rm [\ref{21}].} {\it Для \ $\Omega(t)$ из (\ref{1.5})
с $r>0$ при $0<p<\infty$ верно соотношение}
 $$
\sum\limits_{s\in\chi^{\perp}(N)}\big(\Omega(2^{-s})\big)^{p}\ll
\sum\limits_{s\in\Theta(N)}\big(\Omega(2^{-s})\big)^{p}.
 $$

{\bf Лемма Д \rm [\ref{21}].} {\it Для $\Omega(t)$ из (\ref{1.5}) с
 $0<\beta<r$ при $0<p<\infty$ верно соотношение
 $$
\sum\limits_{s\in\chi^{\perp}(N)}\big(\Omega(2^{-s})2^{\|s\|_{1}\beta}\big)^{p}\ll
\sum\limits_{s\in\Theta(N)}\big(\Omega(2^{-s})2^{\|s\|_{1}\beta}\big)^{p},
 $$
где $\|s\|_{1}=s_1+...+s_d, \ s_j\in \mathbb{N}$.}

Далее положим ${M=|Q(N)|}$, тогда с учетом леммы Б, получим
$$
\begin{array}{cc} \displaystyle M\asymp N^{\frac{1}{r}}\big(\log
N\big)^{-\frac{b_{1}}{r}-...-\frac{b_{d}}{r}+d-1},  \\ \displaystyle
\log M\asymp\log N, \ \ N\asymp M^{r}\big(\log
M\big)^{b_{1}+...+b_{d}-(d-1)r}.
\end{array}
$$

\bf{3. Основные результаты.}\rm

Имеет место

{\bf Теорема 3.1}. {\it  Пусть $1\leq q\leq p<\infty, \ p\geq2,  \
1\leq\theta<\infty$,  а $\Omega(t)$ задана формулой (\ref{1.5}).
Тогда при  $r>0$, справедливы порядковые соотношения
\begin{equation}
d^{\perp}_{M}(B^{\Omega}_{p,\theta}, L_{q})\asymp
d^{B}_{M}(B^{\Omega}_{p,\theta}, L_{q})\asymp M^{-r}\big(\log
M\big)^{-b_{1}-...-b_{d}+(d-1)\big(r+
\big(\frac{1}{2}-\frac{1}{\theta}\big)_{+}\big)},
 \label{3.1}
 \end{equation}
где $a_{+}=\max\{a,0\}$.}

Д\,о\,к\,а\,з\,а\,т\,е\,л\,ь\,с\,т\,в\,о. Сначала получим в
(\ref{3.1}) оценку сверху для поперечника
$d^{\perp}_{M}(B^{\Omega}_{p,\theta}, L_{q})$. С этой целью
рассмотрим приближение функций $f\in B^{\Omega}_{p,\theta}$
тригонометрическими полиномами $t_{Q(N)}$ вида
 $$
t_{Q(N)}(x)=\sum\limits_{s\in\chi(N)}\delta_{s}(f,x).
 $$

Воспользовавшись следствием из теоремы А, получим
 $$
\|f(\cdot)-t_{Q(N)}(\cdot)\|_{q}\leq\|f(\cdot)-t_{Q(N)}(\cdot)\|_{p}\ll\bigg(\sum\limits_{s\in
\
\chi^{\perp}(N)}\|\delta_{s}(f,\cdot)\|_{p}^{2}\bigg)^{\frac{1}{2}}=
 $$
 $$
 =\bigg(\sum\limits_{s\in \
\chi^{\perp}(N)}\Omega^{-2}(2^{-s})\|\delta_{s}(f,\cdot)\|_{p}^{2} \
\Omega^{2}(2^{-s})\bigg)^{\frac{1}{2}}=I_{1}.
 $$

Чтобы получить оценку $I_{1}$ рассмотрим два случая.

Пусть сначала  $1\leq\theta\leq2$. Тогда в силу неравенства
[\ref{22}]
 $$
\bigg(\sum\limits_{k}|a_{k}|^{\nu_{2}}\bigg)^{\frac{1}{\nu_{2}}}\leq
\bigg(\sum\limits_{k}|a_{k}|^{\nu_{1}}\bigg)^{\frac{1}{\nu_{1}}},
 \ \ \ 1\leq\nu_{1}\leq\nu_{2}<\infty,
 $$
имеем
 $$
I_{1}\leq\bigg(\sum\limits_{s\in \
\chi^{\perp}(N)}\Omega^{-\theta}(2^{-s})\|\delta_{s}(f,\cdot)\|_{p}^{\theta}
\ \Omega^{\theta}(2^{-s})\bigg)^{\frac{1}{\theta}}\leq
 $$
 $$
\leq N^{-1}\bigg(\sum\limits_{s\in \
\chi^{\perp}(N)}\Omega^{-\theta}(2^{-s})\|\delta_{s}(f,\cdot)\|_{p}^{\theta}
\bigg)^{\frac{1}{\theta}}\ll
 $$
 $$
\ll N^{-1} \|f\|_{B^{\Omega}_{p,\theta}}\asymp M^{-r}\big(\log
M\big)^{-b_{1}-...-b_{d}+(d-1)r}.
 $$

Если же $2<\theta<\infty$, то применив к $I_{1}$ неравентсво
Гельдера с показателем $\frac{\theta}{2}$, получим
 $$
I_{1}\leq\bigg(\sum\limits_{s\in \
\chi^{\perp}(N)}\Omega^{-\theta}(2^{-s})
\|\delta_{s}(f,\cdot)\|_{p}^{\theta}
\bigg)^{\frac{1}{\theta}}\bigg(\sum\limits_{s\in \
\chi^{\perp}(N)}\big(\Omega(2^{-s})\big)^{\frac{2\theta}{\theta-2}}
\bigg)^{\frac{1}{2}-\frac{1}{\theta}}\ll
 $$
 $$
\ll\|f\|_{B_{p,\theta}^{\Omega}}\bigg(\sum\limits_{s\in \
\chi^{\perp}(N)}\big(\Omega(2^{-s})\big)^{\frac{2\theta}{\theta-2}}
\bigg)^{\frac{1}{2}-\frac{1}{\theta}}\leq \bigg(\sum\limits_{s\in \
\chi^{\perp}(N)}\big(\Omega(2^{-s})\big)^{\frac{2\theta}{\theta-2}}
\bigg)^{\frac{1}{2}-\frac{1}{\theta}}=I_{2}.
 $$

Далее, воспользовавшись сначала леммой Г, а затем леммой В,
продолжим оценку~$I_{2}$
 $$
I_{2}\ll\bigg(\sum\limits_{s\in
\Theta(N)}\big(\Omega(2^{-s})\big)^{\frac{2\theta}{\theta-2}}
\bigg)^{\frac{1}{2}-\frac{1}{\theta}}\ll
N^{-1}\bigg(\sum\limits_{s\in
\Theta(N)}1\bigg)^{\frac{1}{2}-\frac{1}{\theta}}\asymp
 $$
 $$
\asymp M^{-r}\big(\log M\big)^{-b_{1}-...-b_{d}+(d-1)r}\big(\log
M\big)^{(d-1)\big(\frac{1}{2}-\frac{1}{\theta}\big)}=
 $$
 $$
= M^{-r}\big(\log
M\big)^{-b_{1}-...-b_{d}+(d-1)\big(r+\frac{1}{2}-\frac{1}{\theta}\big)}.
 $$

Оценка сверху в (\ref{3.1}) установлена.

Перейдем к доказательству соответствующей оценки снизу, при этом мы
будем придерживаться схемы, примененной В.Н. Темляковым в работе
[\ref{10}]. Заметим, что поскольку полученная оценка сверху не
зависит от параметра $q$, то для доказательства оценки снизу
величины $d_M^{B}(B^{\Omega}_{p,\theta},L_q)$ достаточно рассмотреть
случай $q=1$. Доказательство разобьем на две части.

Пусть сначала $1\leq\theta<2$. В этом случае будем использовать
пример 1 из [\ref{10}]. \ Итак, \ пусть \ число $M$ задано, $G\in
\mathcal{L}_{M}(B)_{1}$. Тогда существует вектор
$k^{0}=(k_{1}^{0},..., k^{0}_{d})\in \widetilde{Q}(N)$, где
$\widetilde{Q}(N)=\bigcup\limits_{s\in\Theta(N)}\rho(s)$ такой, что

\begin{equation}
\big\|e^{i(k^{0},\cdot)}-Ge^{i(k^{0},\cdot)}\big\|_{1}\gg1.
 \label{3.2}
 \end{equation}

Теперь рассмотрим функцию
 $$
g_{1}(x)=N^{-1}e^{i(k^{0},x)},
 $$
которая, как легко видеть, принадлежит классу
$B_{p,\theta}^{\Omega}, \ 1\leq\theta<2$.

Далее, воспользовавшись соотношением (\ref{3.2}), получим
 $$
\big\|g_{1}(\cdot)-Gg_{1}(\cdot)\big\|_{1}=N^{-1}
\big\|e^{i(k^{0},\cdot)}-Ge^{i(k^{0},\cdot)}\big\|_{1}\gg
 $$
 $$
\gg N^{-1}\asymp M^{-r}\big(\log M\big)^{-b_{1}-...-b_{d}+(d-1)r}.
 $$

Для установления оценки снизу величины
$d^{B}_{M}(B^{\Omega}_{p,\theta}, L_{1})$ в случае
$2\leq\theta<\infty$, рассмотрим функцию, аналогичную функции из
примера 6 работы [\ref{10}].

С помощью тех же рассуждений, что и в [\ref{23}], можно показать,
что существует множество $\Theta'(N)\subset\Theta(N)$ такое, что для
$s=(s_1,...,s_d)\in\Theta'(N)$ будет
 $$
 s_j\asymp\log N, \ \ j=\overline{1,d} \ \ \mbox{и} \ \ |\Theta'(N)|\asymp \big(\log N\big)^{d-1}.
 $$

Далее, для $G\in \mathcal{L}_{M}(B)_{1}$ найдутся $N, \
 \Theta_{1}(N)\subset\Theta'(N)$ такие, что
 $$
|\Theta_{1}(N)|\geq\frac{1}{2}|\Theta'(N)|,
 $$
и в каждом $\rho(s), \ s\in\Theta_{1}(N)$, найдутся векторы $k^{s}$
такие, что для функции
 $$
g_{2}(x)=\sum\limits_{s\in\Theta_{1}(N)}e^{i(k^{s},x)}
 $$
 найдется $y^{*}=(y^{*}_{1},..., y^{*}_{d})$ такой, что
 \begin{equation}
\|g_{2}(x+y^{*})-Gg_{2}(x+y^{*})\|_{1}\gg\big(\log
M\big)^{\frac{d-1}{2}}.
 \label{3.3}
 \end{equation}

Доказательство соотношения (\ref{3.3}) проводится по той же схеме,
что и доказательство примера 6 из [\ref{10}].

Итак, рассмотрим функцию
 $$
 g_{3}(x)=C_{5}N^{-1}\big(\log
N\big)^{-\frac{d-1}{\theta}}g_{2}(x), \ C_{5}>0.
 $$

 Легко убедиться, что при соответствующем выборе постоянной $C_{5}$
 функция $g_3$
 принадлежит классу $B_{p,\theta}^{\Omega}, 2\leq\theta<\infty$.

Действительно,
 $$
\|g_{3}\|_{B_{p,\theta}^{\Omega}}=\bigg(\sum
\limits_{s\in\Theta_{1}(N)}
\Omega^{-\theta}(2^{-s})\|\delta_{s}(g_{3},\cdot)\|_{p}^{\theta}
\bigg)^{\frac{1}{\theta}}\ll
 $$
 $$
\ll N^{-1}\big(\log N\big)^{-\frac{d-1}{\theta}} \bigg(\sum
\limits_{s\in\Theta_{1}(N)}
\Omega^{-\theta}(2^{-s})\|\delta_{s}(g_{2},\cdot)\|_{p}^{\theta}
\bigg)^{\frac{1}{\theta}}\ll
 $$
 $$
\ll N^{-1}\big(\log N\big)^{-\frac{d-1}{\theta}}\bigg(\sum
\limits_{s\in\Theta_{1}(N)}
\Omega^{-\theta}(2^{-s})\bigg)^{\frac{1}{\theta}}\ll
  N^{-1}\big(\log
 N\big)^{-\frac{d-1}{\theta}}\cdot N\cdot
 |\Theta_{1}(N)|^{\frac{1}{\theta}}
  \asymp
 $$
 $$
  \asymp\big(\log
 N\big)^{-\frac{d-1}{\theta}}\big(\log
 N\big)^{\frac{d-1}{\theta}}=1.
 $$

Таким образом, воспользовавшись соотношением (\ref{3.3}), будем
иметь
 $$
\|g_{3}(x+y^{*})-Gg_{3}(x+y^{*})\|_{q}\geq
\|g_{3}(x+y^{*})-Gg_{3}(x+y^{*})\|_{1} \gg
 $$
 $$
 \gg N^{-1}\big(\log N\big)^{-\frac{d-1}{\theta}} \|g_{2}(x+y^{*})-Gg_{2}(x+y^{*})\|_{1}\gg
 $$
 $$
 \gg N^{-1}\big(\log
 N\big)^{-\frac{d-1}{\theta}}\big(\log M\big)^{\frac{d-1}{2}}\asymp M^{-r}\big(\log
M\big)^{-b_{1}-...-b_{d}+(d-1)\big(r+\frac{1}{2}-\frac{1}{\theta}\big)}.
 $$

Оценки снизу величины  $d^{B}_{M}(B^{\Omega}_{p,\theta}, L_{1})$ в
обоих случаях установлены. Теорема доказана.

{\bf Теорема 3.2}. {\it Пусть $1\leq q\leq p\leq2, \ (p,q)\neq(1,1),
\ 1\leq\theta<\infty$, а $\Omega(t)$ задана формулой (\ref{1.5}).
Тогда при $r>0$ справедливы оценки}
 \begin{equation}
d^{\perp}_{M}(B^{\Omega}_{p,\theta}, L_{q})\asymp
d^{B}_{M}(B^{\Omega}_{p,\theta}, L_{q})\asymp M^{-r}\big(\log
M\big)^{-b_{1}-...-b_{d}+(d-1)\big(r+\big(\frac{1}{p}-\frac{1}{\theta}\big)_{+}\big)}.
 \label{3.4}
 \end{equation}

 Д\,о\,к\,а\,з\,а\,т\,е\,л\,ь\,с\,т\,в\,о.
Оценка сверху в соотношении (\ref{3.4}) следует из рассуждений
аналогичных тем, которые проводились при получении оценки сверху в
теореме 3.1.

Переходя к установлению оценки снизу для
$d^{B}_{M}(B^{\Omega}_{p,\theta}, L_{q})$, заметим, что ее
достаточно установить для случая $q=1, \ 1<p\leq2.$ Кроме того, если
${\theta\in[1,p)}$, то оценка снизу величины
$d^{B}_{M}(B^{\Omega}_{p,\theta}, L_{1})$ устанавливается с помощью
тех же рассуждений, которые проводились при доказательстве оценки
снизу в теореме~3.1 в случае $1\leq\theta<2$. Таким образом,
остановимся на рассмотрении случая $p\leq\theta<\infty$.

Рассмотрим функцию, аналогичную функции из примера 7 работы
[\ref{10}].

Пусть $N$ достаточно велико. Положим
 $$
 v=[|\Theta'(N)|^{\frac{1}{d}}],
 $$
где множество $\Theta'(N)$ определено выше, $[a]$ --- целая часть
числа $a$.

Далее, разобьем куб $\pi_{d}$ на $v^{d}$ кубов с длиной ребра
$\frac{2\pi}{v}$. Установим взаимно однозначное соответствие между
множеством $\overline{\Theta}(N)\subset\Theta'(N),
|\overline{\Theta}(N)|=v^{d}$ и получившимся множеством кубов. При
этом через $x^{s}\in\pi_{d}$ обозначим центр куба, соответствующего
вектору  $s\in\overline{\Theta}(N)$, и положим
 $$
u=2^{\big[\frac{1}{d}\log|\Theta'(N)|\big]}\asymp\big(\log
N\big)^{\frac{d-1}{d}}.
 $$

Пусть далее $K_n(t)$ обозначает ядро Фейера порядка $n$, т.е.
 $$
K_n(t)=1+2\sum_{k=1}^n\Big(1-\frac{k}{n+1}\Big)\cos kt.
 $$

Через $k^s$ обозначим вектор $k^s=(k_1^{s_1},..., k_d^{s_d})$, где
 $$
k_j^{s_j}= \left\{\begin{array}{cc} \displaystyle
2^{s_j-1}+2^{s_j-2},\ & s_j\geq 2;\\
1,\ & s_j=1,j=\overline{1,d}.
\end{array} \right.
 $$

Пусть $G\in \mathcal{L}_{M}(B)_{1}$. Тогда существуют число $N$ и
множество $\Theta_{2}(N)\subset\overline{\Theta}(N)$ такие, что
$$
|\Theta_{2}(N)|\geq\frac{1}{2}|\overline{\Theta}(N)|,
$$
 и в каждом $\rho(s)$, ${s\in\Theta_{2}(N)}$, найдутся кубы с центром в
$k^{s}$ и длинами ребер $2u$ такие, что для функции
 $$
g_{4}(x)=\sum_{s\in \Theta_{2}(N)} e^{i(k^s,x)}\prod_{j=1}^d
K_u(x_j-x_j^s)
 $$
и некоторого вектора $y^*$ имеет место оценка
 \begin{equation}
 \left\|g_{4}(x+y^*)-Gg_{4}(x+y^*)\right\|_1\gg\big(\log M\big)^{d-1}.
 \label{3.5}
 \end{equation}

Доказательство оценки (\ref{3.5}) проводится с помощью тех же
рассуждений, которые проводились при доказательстве соответствующей
оценки в примере 7 из работы~[\ref{10}].

Теперь возвратимся непосредственно к оценке снизу величины
$d^{B}_{M}(B^{\Omega}_{p,\theta}, L_{1})$. Рассмотрим функцию
 $$
 g_{5}(x)=C_{6}N^{-1}\big(\log
 N\big)^{(d-1)\big(\frac{1}{p}-1-\frac{1}{\theta}\big)} g_{4}(x), \ \ \ C_{6}>0,
 $$
и оценим $||g_{5}||_{B^{\Omega}_{p,\theta}}$.

Учитывая, что в силу выбора параметра $u$
 $$
\left\|A_s(g_{4},\cdot)\right\|_p\ll\big\|K_u(\cdot) \big\|_p \asymp
u^{d\big(1-\frac{1}{p}\big)}\asymp \big(\log
 N\big)^{(d-1)\big(1-\frac{1}{p}\big)}, \ \ \ \forall s\in \Theta_{2}(N),
 $$
будем иметь
 $$
||g_{5}||_{B^{\Omega}_{p,\theta}}=\bigg(\sum_{s\in
\Theta_{2}(N)}\Omega^{-\theta}(2^{-s})||A_s(g_{5},\cdot)||^{\theta}_p\bigg)
^{\frac{1}{\theta}}\ll
 $$
 $$
\ll N^{-1} \big(\log
N\big)^{(d-1)\big(\frac{1}{p}-1-\frac{1}{\theta}\big)}\cdot
\bigg(\sum_{s\in \Theta_{2}(N)} \Omega^{-\theta}(2^{-s})
||A_s(g_{4},\cdot)||^{\theta}_p \bigg) ^{\frac{1}{\theta}}\ll
 $$
 $$
\ll N^{-1} \big(\log
N\big)^{(d-1)\big(\frac{1}{p}-1-\frac{1}{\theta}\big)}
\cdot\big(\log N\big)^{(d-1)\big(1-\frac{1}{p}\big)}
\cdot\bigg(\sum_{s\in \Theta_{2}(N)} \Omega^{-\theta}(2^{-s})\bigg)
^{\frac{1}{\theta}}\asymp
 $$
\begin{equation}
\asymp \big(\log N\big)^{-\frac{d-1}{\theta}}\cdot
|\Theta_{2}(N)|^{\frac{1}{\theta}}\asymp1. \label{3.6}
\end{equation}

Таким образом, из (\ref{3.6}) заключаем, что функция $g_{5}\in
B^{\Omega}_{p,\theta}, \ p\leq\theta<\infty,$ с соответствующей
постоянной $C_{6}>0$.

Далее, согласно (\ref{3.5}), существует вектор $y^*$  такой, что
 $$
\left\|g_{5}(x+y^*)-Gg_{5}(x+y^*)\right\|_1\gg
 $$
 $$\gg N^{-1} \big(\log N\big)^{(d-1)\big(\frac{1}{p}-1-\frac{1}{\theta}\big)} \left\|g_{4}(x+y^*)-Gg_{4}(x+y^*)\right\|_1\gg
 $$
 $$
\gg N^{-1} \big(\log
N\big)^{(d-1)\big(\frac{1}{p}-\frac{1}{\theta}\big)}\asymp
M^{-r}\big(\log
M\big)^{-b_{1}-...-b_{d}+(d-1)\big(r+\frac{1}{p}-\frac{1}{\theta}\big)}.
 $$
Теорема доказана.

{\bf Теорема 3.3}. {\it Пусть $1\leq p<\infty,  \
1\leq\theta<\infty$,  функция $\Omega(t)$ задана формулой
(\ref{1.5}). Тогда при $r>\frac{1}{p}$ справедливы соотношения}
 \begin{equation}
d^{\perp}_{M}(B^{\Omega}_{p,\theta}, L_{\infty})\asymp
d^{B}_{M}(B^{\Omega}_{p,\theta}, L_{\infty})\asymp
  M^{-r+\frac{1}{p}} \big(\log M\big)^{-b_{1}-...-b_{d}+(d-1)\big(r+1-\frac{1}{p}-\frac{1}{\theta}\big)}.
 \label{3.7}
 \end{equation}

Д\,о\,к\,а\,з\,а\,т\,е\,л\,ь\,с\,т\,в\,о. Сначала установим в
(\ref{3.7}) оценку сверху. Для этого достаточно получить
соответствующую оценку сверху для величины
$\|f(\cdot)-t_{Q(N)}(\cdot)\|_{\infty}$, где $f\in
B^{\Omega}_{p,\theta}$, \
$t_{Q(N)}(x)=\sum\limits_{s\in\chi(N)}\delta_{s}(f,x)$.

Итак, пусть $q_{0}$ --- некоторое число, удовлетворяющее условию
$1<q_{0}<\infty$. Тогда, воспользовавшись неравенством Минковского и
затем неравенством разных метрик Никольского, для $f\in
B^{\Omega}_{p,\theta}$ будем иметь
$$
\|f(\cdot)-t_{Q(N)}(\cdot)\|_{\infty}=
\bigg\|f(\cdot)-\sum\limits_{s\in
 \chi(N)}\delta_{s}(f,\cdot)\bigg\|_{\infty}\leq
\sum\limits_{s\in
\chi^{\perp}(N)}\|\delta_{s}(f,\cdot)\|_{\infty}\ll
 $$
 $$
 \ll \sum\limits_{s\in
\chi^{\perp}(N)}2^{\frac{\|s\|_{1}}{q_{0}}}\|\delta_{s}(f,\cdot)\|_{q_{0}}
\asymp \sum\limits_{s\in
\chi^{\perp}(N)}2^{\frac{\|s\|_{1}}{q_{0}}}\|A_{s}(f,\cdot)\|_{q_{0}}
\ll
 $$
 $$
 \ll \sum\limits_{s\in
\chi^{\perp}(N)}2^{\frac{\|s\|_{1}}{q_{0}}}2^{\|s\|_{1}\big(\frac{1}{p}
-\frac{1}{q_0}\big)} \|A_{s}(f,\cdot)\|_{p} =
 $$
 $$
  =\sum\limits_{s\in
\chi^{\perp}(N)}\Omega(2^{-s})2^{\frac{\|s\|_{1}}{p}}\Omega^{-1}(2^{-s})
\|A_{s}(f,\cdot)\|_{p}=I_{3}.
 $$

Применив к последней сумме в $I_{3}$ неравенство Гельдера с
показателем  $\theta$ (с естественной модификацией при $\theta=1$) и
воспользовавшись леммой Д, продолжим оценку $I_{3}$
 $$
I_{3}\leq
\bigg(\sum\limits_{s\in\chi^{\perp}(N)}\Omega^{-\theta}(2^{-s})
\|A_{s}(f,\cdot)\|_{p}^{\theta}\bigg)^{\frac{1}{\theta}}\bigg(\sum\limits_{s\in
\chi^{\perp}(N)}\big(\Omega(2^{-s})2^{\frac{\|s\|_{1}}{p}}\big)^{\frac{\theta}{\theta-1}}
\bigg)^{1-\frac{1}{\theta}}\ll
 $$
 $$
\ll\|f\|_{B^{\Omega}_{p,\theta}}\bigg(\sum\limits_{s\in
\chi^{\perp}(N)}\big(\Omega(2^{-s})2^{\frac{\|s\|_{1}}{p}}\big)^{\frac{\theta}{\theta-1}}
\bigg)^{1-\frac{1}{\theta}}\ll
 $$
 $$
\ll\bigg(\sum\limits_{s\in
\Theta(N)}\big(\Omega(2^{-s})2^{\frac{\|s\|_{1}}{p}}\big)^{\frac{\theta}{\theta-1}}
\bigg)^{1-\frac{1}{\theta}}=I_{4}.
 $$

Далее, учитывая тот факт, что для $s\in\Theta(N)$ будет
 $$
2^{\|s\|_1}\asymp
N^{\frac{1}{r}}\prod\limits_{j=1}^{d}s_{j}^{-\frac{b_j}{r}},
 $$
продолжим оценку $I_4$
 $$
I_4\ll N^{-1}\cdot N^{\frac{1}{pr}}\cdot\bigg(\sum\limits_{s\in
\Theta(N)}\bigg(\prod\limits_{j=1}^{d}s_{j}^{-\frac{b_j}{pr}}\bigg)^
{\frac{\theta}{\theta-1}}\bigg)^{1-\frac{1}{\theta}}\ll
 $$
 $$
\ll N^{\frac{1}{pr}-1}\cdot \big(\log N\big)^{-\frac{b_1}{pr}-
...-\frac{b_d}{pr}}\bigg(\sum\limits_{s\in
\Theta(N)}1\bigg)^{1-\frac{1}{\theta}}\asymp
 $$
 $$
\asymp \Big(M^{r}\big(\log
M\big)^{b_{1}+...+b_{d}-(d-1)r}\Big)^{\frac{1}{pr}-1}\cdot \big(\log
M\big)^{-\frac{b_1}{pr}- ...-\frac{b_d}{pr}}\cdot \big(\log
M\big)^{(d-1)\big(1-\frac{1}{\theta}\big)}=
 $$
 $$
 =M^{-r+\frac{1}{p}} \big(\log M\big)^{-b_{1}-...-b_{d}+(d-1)\big(r+1-\frac{1}{p}-\frac{1}{\theta}\big)}.
 $$

Таким образом, оценка сверху в (\ref{3.7}) установлена.

Переходя к оценке снизу, напомним, что ее достаточно получить для
величини $d^{B}_{M}(B^{\Omega}_{p,\theta}, L_{\infty})$, при этом мы
будем придерживаться схемы, примененной в примере 2 из [\ref{10}].

Рассмотрим функцию
$$
g_6(x)=\sum\limits_{s\in \Theta'(N)}\mathcal{K}^{s}(x),
$$
где множество $\Theta'(N)$ --- определено выше, а
 $$
\mathcal{K}^{s}(x)=e^{i(k^{s},x)}
\prod\limits_{j=1}^{d}K_{2^{s_{j}-2}}(x_{j}),
 $$
$K_n(t)$ --- ядро Фейера.

   Предположим, что оператор $G$ принадлежит $\mathcal{L}_M(B)_\infty$. Покажем, что существует вектор $y^*=(y_1^*,..., y_d^*)$ такой, что
\begin{equation}
\|g_6(x-y^*) - Gg_6(x-y^*)\|_\infty\gg M. \label{3.8}
\end{equation}

Очевидно,
$$
\|g_6(x-y) - Gg_6(x-y)\|_\infty\geq g_6(0) - \min\limits_{y=x}
\mbox{Re} \ Gg_6(x-y).
$$

Воспользовавшись леммой А, будем иметь
\begin{equation}
\min\limits_{y=x} \mbox{Re} \ Gg_6(x-y)\leq
M^{\frac{1}{2}}B\bigg(\sum\limits_k|\widehat{g}_6(k)|^2
\bigg)^{\frac{1}{2}}\ll
M^{\frac{1}{2}}B|\widetilde{Q}'(N)|^{\frac{1}{2}}. \label{3.9}
\end{equation}

Здесь через $|\widetilde{Q}'(N)|$ обозначено количество элементов
множества
$$
\widetilde{Q}'(N)=\bigcup\limits_{s\in \Theta'(N)}\rho(s).
$$

Далее, учитывая, что $|\Theta'(N)|\asymp\big(\log N\big)^{d-1}$, а
также соотношение
$$
|\rho(s)|=2^{\|s\|_1}\asymp N^{\frac{1}{r}}\big(\log
N\big)^{-\frac{b_{1}}{r}-...-\frac{b_{d}}{r}},
$$
будем иметь
\begin{equation}
|\widetilde{Q}'(N)|\asymp N^{\frac{1}{r}}\big(\log
N\big)^{-\frac{b_{1}}{r}-...-\frac{b_{d}}{r}+d-1}. \label{3.10}
\end{equation}

С другой стороны,
\begin{equation}
g_6(0)\asymp N^{\frac{1}{r}}\big(\log
N\big)^{-\frac{b_{1}}{r}-...-\frac{b_{d}}{r}+d-1}\asymp
|\widetilde{Q}'(N)|. \label{3.11}
\end{equation}

С учетом (\ref{3.9}) и (\ref{3.10}) можно подобрать такое $N$, что
${|\widetilde{Q}'(N)|\asymp M}$ и правая часть в (\ref{3.11}) будет
по крайней мере вдвое больше правой части (\ref{3.9}). Для этого $N$
при некотором $y^*$, будем иметь
$$
\|g_6(x-y^*) - Gg_6(x-y^*)\|_\infty\gg M.
$$

Таким образом оценка (\ref{3.8}) доказана.

Рассмотрим функцию
 $$
 g_{7}(x)=C_{7}N^{-1} \Big(N^{\frac{1}{r}}\big(\log
 N\big)^{-\frac{b_{1}}{r}-...-\frac{b_{d}}{r}}\Big)^{\frac{1}{p}-1} \big(\log N\big)^{-\frac{d-1}{\theta}}g_6(x), \
 C_{7}>0.
 $$

Покажем, что функция $g_{7}$ при надлежащем выборе постоянной $C_7$
принадлежит классу $B^{\Omega}_{p,\theta}, \ 1\leq p<\infty,  \
1\leq\theta<\infty$.

Действительно, поскольку в силу свойств ядра Фейера
 $$
\left\|\mathcal{K}^s\right\|_p\asymp
2^{\|s\|_1\big(1-\frac{1}{p}\big)}, \ 1\leq p\leq\infty,
 $$
то будем иметь
 $$
||g_{7}||_{B^{\Omega}_{p,\theta}}=\bigg(\sum_{s\in
\Theta'(N)}\Omega^{-\theta}(2^{-s})||A_s(g_{7},\cdot)||^{\theta}_p\bigg)
^{\frac{1}{\theta}}\ll
 $$
 $$
\ll N^{-1}\Big(N^{\frac{1}{r}}\big(\log
 N\big)^{-\frac{b_{1}}{r}-...-\frac{b_{d}}{r}}\Big)^{\frac{1}{p}-1} \big(\log N\big)^{-\frac{d-1}{\theta}}\bigg(\sum_{s\in \Theta'(N)} \Omega^{-\theta}(2^{-s})
||A_s(g_6,\cdot)||^{\theta}_p \bigg) ^{\frac{1}{\theta}}\asymp
 $$
 $$
\asymp \Big(N^{\frac{1}{r}}\big(\log
 N\big)^{-\frac{b_{1}}{r}-...-\frac{b_{d}}{r}}\Big)^{\frac{1}{p}-1} \big(\log N\big)^{-\frac{d-1}{\theta}}
 \bigg(\sum_{s\in \Theta'(N)} 2^{\|s\|_1\big(1-\frac{1}{p}\big)\theta}\bigg)
^{\frac{1}{\theta}}\asymp
 $$
 $$
\asymp \Big(N^{\frac{1}{r}}\big(\log
 N\big)^{-\frac{b_{1}}{r}-...-\frac{b_{d}}{r}}\Big)^{\frac{1}{p}-1} \big(\log N\big)^{-\frac{d-1}{\theta}}\Big(N^{\frac{1}{r}}\big(\log
 N\big)^{-\frac{b_{1}}{r}-...-\frac{b_{d}}{r}}\Big)^{1-\frac{1}{p}} |\Theta'(N)|^{\frac{1}{\theta}}\asymp
 $$
 $$
\asymp \big(\log N\big)^{-\frac{d-1}{\theta}}\big(\log
N\big)^{\frac{d-1}{\theta}}=1.
 $$

Теперь, воспользовавшись соотношением (\ref{3.8}), получим
 $$
\left\|g_{7}(x-y^*)-Gg_{7}(x-y^*)\right\|_\infty\gg
 $$
 $$
 \gg N^{-1} \Big(N^{\frac{1}{r}}\big(\log
 N\big)^{-\frac{b_{1}}{r}-...-\frac{b_{d}}{r}+d-1}\Big)^{\frac{1}{p}-1} \big(\log N\big)^{(d-1)\big(1-\frac{1}{p}-\frac{1}{\theta}\big)} \times
 $$
 $$
  \times \left\|g_6(x-y^*)-Gg_6(x-y^*)\right\|_\infty\gg
 $$
 $$
  \gg M^{-r}\big(\log
M\big)^{-b_{1}-...-b_{d}+(d-1)r}M^{\frac{1}{p}-1}\big(\log
M\big)^{(d-1)\big(1-\frac{1}{p}-\frac{1}{\theta}\big)}M=
 $$
  $$
 = M^{-r+\frac{1}{p}} \big(\log M\big)^{-b_{1}-...-b_{d}+(d-1)\big(r+1-\frac{1}{p}-\frac{1}{\theta}\big)}.
 $$

Оценка снизу для $d^{B}_{M}(B^{\Omega}_{p,\theta}, L_{\infty})$, а
также и для поперечника $d^{\perp}_{M}(B^{\Omega}_{p,\theta},
L_{\infty})$ установлена. Теорема доказана.

Подытоживая полученные результаты отметим, что оптимальными (в
смысле порядка) подпространствами в теоремах 3.1 -- 3.3 являются
подпространства тригонометрических полиномов с "номерами" \ гармоник
из множеств $Q(N)$.

{\bf Замечание 3.1.} {\it При $d=2$ и соответствующих ограничениях
на  параметры $p$ и $q$, результаты теорем 3.1 -- 3.3 получены в
[\ref{1}, \ref{2}].}

{\bf Замечание 3.2.} {\it В случае, когда
$\Omega(t)=\prod\limits_{j=1}^dt_j^{r},$ результаты теорем 3.1 --
3.3 (для классов $B^{r}_{p,\theta}, 1\leq\theta<\infty$) получены
А.С. Романюком в [\ref{15a}, \ref{15}].}

{\bf Замечание 3.3.} {\it Порядки поперечников
$d^{\perp}_{M}(H^{\Omega}_{p}, L_{q})$ и величин
$d^{B}_{M}(H^{\Omega}_{p}, L_{q})$ при $p$ и $q$, которые
удовлетворяют условиям теорем 3.1 -- 3.3 получены Н.Н.~Пустовойтовым
в [\ref{14}].}

Выражаю искреннюю благодарность А.С. Романюку за полезные замечания
и ценные советы, сделанные им в процессе обсуждения результатов
работы.

\newpage
\begin{enumerate}
\Rus

\item \label{1} {\it Конограй А.\,Ф.\/} Оцінки апроксимативних
характеристик класів $B^{\Omega}_{p,\theta}$ періодичних функцій
двох змінних з заданою мажорантою мішаних модулів неперервності //
Укр. мат. журн. --- 2011. --- {\bf 63}, №2. --- C.~176--186.

\item \label{2} {\it Конограй А.\,Ф.\/} Оцінки апроксимативних
характеристик класів $B^{\Omega}_{p,\theta}$ періодичних функцій
двох змінних з заданою мажорантою мішаних модулів неперервності в
просторі $L_\infty$ // Теорія наближення функцій та суміжні питання:
Зб. праць Ін-ту математики НАН України. --- 2011. --- {\bf 8}, №1.
--- C.~97--110.

\item \label{3} {\it Бари Н.К., Стечкин С.Б.} Наилучшие
приближения и дифференциальные свойства двух сопряженных функций
// Тр. Моск. мат. о-ва. --- 1956. --- {\bf 5}.~--- C.~483 -- 522.

\item \label{4} {\it Sun Youngsheng, Wang Heping. \/}
Representation and approximation of multivariate periodic functions
with bounded mixed moduli of smoothness // Тр. Мат. ин-та им.
В.А.~Стеклова.
--- 1997. --- {\bf 219}. --- С.~356 -- 377.

\item \label{5} {\it Пустовойтов
Н.Н.\/} Представление и приближение периодических функций многих
переменных с заданным смешанным модулем непрерывности
// Anal. Math. --- 1994. --- {\bf 20}, № 1. --- P.~35 -- 48.

\item \label{6} {\it Лизоркин
П.И., Никольский С.М.\/} Пространства функций смешанной гладкости с
декомпозиционной точки зрения // Тр. Мат. ин-та им. В.А. Стеклова АН
СССР. --- 1989.--- {\bf 187}.--- C.~143--161.

\item \label{7}{\it Стасюк С.А., Федуник О.В.\/} Апроксимативні
характеристики класів $B_{p,\theta}^{\Omega}$ періодичних функцій
багатьох змінних
// Укр. мат. журн. --- 2006. --- {\bf 58}, №~5. --- C.~692 -- 704.

\item \label{8}{\it  Темляков В.\,Н.\/} Поперечники некоторых
классов функций нескольких переменных // ДАН СССР. --- 1982. ---
{\bf 267}, №~2. --- C.~314--317.

\item \label{9}{\it Динь Зунг \/}  Приближение функций многих
переменных на торе тригонометрическими полиномами //  Матем. сб. ---
1986. --- {\bf 131}, №~2. --- C.~251--271.

\item \label{10}{\it Темляков В.\,Н.\/} Оценки асимптотических
характеристик классов функций с ограниченной смешанной производной
или разностью // Тр. МИАН СССР. --- 1989. ---  {\bf 189}. ---
C.~138--168.

\item \label{11}{\it Галеев Э.\,М.\/} Порядки ортопроекционных
поперечников классов периодических функций одной и нескольких
переменных // Матем. заметки. --- 1988. --- {\bf 43}, №~2. ---
C.~197--211.

\item \label{12}{\it Галеев Э.\,М.\/} Приближение классов
периодических функций нескольких переменных ядерными операторами //
 Матем. заметки. ---1990. --- {\bf 47}, №~3. --- C.~32--41.

\item \label{13}{\it Андрианов А.\,В., Темляков В.\,Н.\/} О двух
методах распространения свойств систем функций от одной переменной
на их тензорное произведение // Тр. МИРАН. --- 1997. --- {\bf 47}.
--- C.~32--43.

\item \label{14}{\it Пустовойтов Н.\,Н.\/} Ортопоперечники классов
многомерных периодических функций, мажоранта смешанных модулей
непрерывности которых содержит как степенные, так и логарифмические
множители // Anal. Math. --- 2008. --- {\bf 34}. --- C.~187--224.

\item \label{15a}{\it Романюк А.\,С.\/} Наилучшие приближения и
поперечники классов периодических функций многих переменных //
Матем. сб. --- 2008. --- {\bf 199}, №~2. --- C.~93--114.

\item \label{15}{\it Романюк А.\,С.\/} Поперечники и наилучшее
приближение классов $B_{p, \theta}^r$ периодических функций многих
переменных // Anal. Math. --- 2011. --- {\bf 37}. --- C.~181--213.

\item \label{16}{\it Темляков В.\,Н.\/} Приближение функций с
ограниченной смешанной производной// Тр. МИАН СССР. --- 1986. ---
{\bf 178}. --- 112~c.

\item \label{17}{\it Temlyakov V.\,N.\/} Approximation of Periodic
Functions. --- NY.: Comput. Math. Anal. Ser., Nova Science Publ,
1993.

\item \label{18}{\it  Никольский С.\,М.\/} Приближение функций многих
переменных и теоремы вложения. --- М.: Наука, 1969. --- 480~c.

\item \label{19}{\it Никольский С.\,М.\/} Неравенства для целых
функций конечной степени и их применение в теории дифференцируемых
функций многих переменных // Тр. МИАН СССР. --- 1951. --- {\bf 38}.
--- C.~244--278.

\item \label{20}{\it Jakson D.\/} Certain problem of closest
approximation // Bull. Amer. Math. Soc. --- 1933. --- {\bf 39},
№~12. --- C.~889--906.

\item \label{21}{\it Пустовойтов Н.\,Н.\/} Приближение многомерных
функций с заданной мажорантой смешанных модулей непрерывности //
Матем. заметки. --- 1999. --- {\bf 65}, №~1. --- C.~107--117.

\item \label{22}{\it Харди Г., Литтлвуд Д., Полиа Г.\/} Неравенства.
--- М.: Изд-во иностр. лит., 1948.

\item \label{23}{\it Пустовойтов Н.\,Н.\/} О приближении и
характеризации периодических функций многих переменных, имеющих
мажоранту смешанных модулей непрерывности специального вида // Anal.
Math. --- 2003. --- {\bf 29}. --- C.~201--218

\end{enumerate}
\end{document}